\documentclass[12pt]{amsart}

\theoremstyle{plain}
\newtheorem{lemma}{Lemma}[section]

\newtheorem*{mainthm}{Main Theorem}

\newtheorem{propo}[lemma]{Proposition}

\newtheorem{con}[lemma]{Conjecture}

\newtheorem{claim}{Claim}

\theoremstyle{definition}
\newtheorem{defi}[lemma]{Definition}
\newtheorem{defs}[lemma]{Definitions}

\theoremstyle{remark}

\newcommand{\pic}{{\rm Pic}}
\newcommand{\bs}{{\rm Bs}}
\newcommand{\p}{\mathbb{P}}

\newcommand{\f}{\mathbb{F}}

\newcommand{\oc}{{\mathcal O}}

\newcommand{\ls}{{\mathcal L}}

\newcommand{\ms}{{\mathcal M}}

\newcommand{\rs}{{\mathcal R}}

\newcommand{\dlp}{l_{\p}}
\newcommand{\dlf}{l_{\f}}
\newcommand{\dklp}{\hat{l}_{\p}}
\newcommand{\dklf}{\hat{l}_{\f}}

\newcommand{\lp}{\ls_{\p}}
\newcommand{\lf}{\ls_{\f}}
\newcommand{\klp}{\hat{\ls}_{\p}}
\newcommand{\klf}{\hat{\ls}_{\f}}

\newcommand{\vp}{v_{\p}}
\newcommand{\vf}{v_{\f}}
\newcommand{\kvp}{\hat{v}_{\p}}
\newcommand{\kvf}{\hat{v}_{\f}}

\newcommand{\rt}{\rightarrow}

\numberwithin{equation}{section}
\numberwithin{figure}{section}
\numberwithin{table}{section}

\begin{document}

\title[Quasi-homogeneous linear systems on $\p^2$]{Quasi-homogeneous linear systems on $\p^2$ with base points of multiplicity $5$}
\author{Antonio Laface and Luca Ugaglia}
\address{Department of Mathematics, University of Milan, Via Saldini 50,
20100 MILANO}
\email{laface@mat.unimi.it\\
ugaglia@mat.unimi.it}
\begin{abstract}
In this paper we consider linear systems of $\p^2$ with all but one of 
the base 
points of multiplicity $5$. We give an explicit way to evaluate the dimensions
of such systems. 
\end{abstract}

\maketitle

\section*{Introduction} 

We consider $r$ points $p_1,\ldots,p_r$, in general position on $\p^2$ and 
to each $p_i$ we associate a natural number $m_i$ 
called the {\em multiplicity} of the point. 
Let $r_j$ be the number of points of multiplicity $m_j$
and let 
$\ls(d,{m_1}^{r_1},\ldots ,{m_k}^{r_k})$ be the linear system 
of curves of degree $d$ with $r_j$ general base points of
multiplicity at least $m_j$ for $j=1\ldots k$. 
A linear system is called {\em homogeneous} if all the multiplicities are
equal and {\em quasi-homogeneous} if all but one of the multiplicities are
equal.

\begin{defs}
The {\em effective dimension} of the system is defined to be
\[
l(d,{m_1}^{r_1},\ldots ,{m_k}^{r_k})=\dim 
\ls(d,{m_1}^{r_1},\ldots ,{m_k}^{r_k}),
\]
the {\em virtual dimension} to be:
\[
v(d,{m_1}^{r_1},\ldots ,{m_k}^{r_k})=
\frac{d(d+3)}{2}-\sum r_i\frac{m_i(m_i+1)}{2},
\]
and the {\em expected dimension}
\[
e=\max\{v,-1\}.
\]
\end{defs}

It follows immediately that for a given system 
\begin{equation}
\label{remark1}
v\leq e\leq l,
\end{equation}
and the second inequality may be strict since
the conditions imposed by the points may fail to be independent. 
In this case we say that the system is {\em special}.

Let $Z$ be the $0$-dimensional scheme defined by the multiple points $p_i$
and consider the exact sequence of sheaves:
\[
0\rt {\mathcal I}(Z)\rt\oc_{S}\rt\oc_Z\rt 0
\]
where ${\mathcal I}(Z)$ is the ideal sheaf of $Z$. 
Tensoring with $L=\oc_{\p^2}(d)$ and taking cohomology we obtain:
\[
0\rt H^0(L\otimes {\mathcal I}(Z))\rt H^0(L)\rt H^0(L_Z)\rt 
H^1(L\otimes {\mathcal I}(Z))\rt 0.
\]
In this way we see that 
\begin{equation}
\label{equx}
h^1(L\otimes {\mathcal I}(Z))=
l(d,{m_1}^{r_1},\ldots ,{m_k}^{r_k}) - 
v(d,{m_1}^{r_1},\ldots ,{m_k}^{r_k})
\end{equation}

Consider the blow up $S$ of $\p^2$ at $p_1,\ldots ,p_r$ and
denote by $\ls$ the strict transform of the system 
$\ls(d,{m_1}^{r_1},\ldots ,{m_k}^{r_k})$. 
Define the virtual and the expected dimension of $\ls$ as those of 
$\ls(d,{m_1}^{r_1},\ldots ,{m_k}^{r_k})$. 
By Riemann-Roch, the virtual dimension of $\ls$ may be given 
in the following way:
\[
v(\ls)=\frac{\ls^2-\ls\cdot K_S}{2},
\]
where $K_S$ denotes the canonical bundle of $S$.
We recall that $E$ is a $(-1)$-curve on $S$ if $E$ is irreducible and 
$E^2=E\cdot K_S=-1$. If $\ls\cdot E=-t<0$, then $\bs\mid\ls\mid$ contains $E$ 
with multiplicity $t$. Let $\ms=\ls-tE$ be the residual system, then: 
\[
v(\ls)=\frac{\ls^2-\ls\cdot K_S}{2}=
\frac{(\ms+tE)^2-(\ms+tE)\cdot K_S}{2}=v(\ms)+\frac{t-t^2}{2}
\]
So if $t\geq 2$ and $\ls$, or equivalently $\ms$, is not empty, then
$l(\ls)=l(\ms)\geq
v(\ms)>v(\ls)$.

This allows us to recall the following definition given in \cite{cm1}: 
\begin{defi}
A linear system $\ls$ is \emph{$(-1)$-special}
if there are $(-1)$-curves $E_1,\ldots,E_r$
such that $\ls\cdot E_j = -n_j$, with $n_j \geq 1$ for every $j$,
$n_j \geq 2$ for some $j$, and the residual system 
$\ms = \ls - \sum_j n_j E_j$
has non-negative virtual dimension $v(\ms) \geq 0$.
\end{defi}

The following conjecture was formulated for the first time in \cite{hi}

\begin{con}
\label{con:general}
A linear system $\ls$ on $\p^2$ is special if and only if it is $(-1)$-special.
\end{con}

The conjecture is known to be true in the following cases:
\begin{itemize}
\item The number of points $r\leq 9$ \cite{harb1}.
\item The system is homogeneous of multiplicity $\leq 12$ \cite{cm2}.
\item All the multiplicities are $\leq 4$ \cite{mig1}
\item The system is quasi-homogeneous with points of multiplicity $\leq 3$ \cite{cm1}.
\item The system is quasi-homogeneous with points of multiplicity $4$ \cite{anto,seibert}.
\item The number of points is a power of $4$ and the system is homogeneous \cite{ev3}.
\end{itemize}

In this paper we prove the conjecture in the case of quasi-homogeneous
systems of multiplicity $5$. \\
The paper is organized as follows: in Section 1 we find a complete list of
$(-1)$-curves that may produce quasi-homogeneous $(-1)$-special systems of
multiplicity 5. Section 2 is devoted to the classification of such
systems. In Section 3 we recall the degeneration of $\p^2$ presented in
\cite{cm1}. In Section 4 we prove the main theorem under the assumption of
some technical lemmas, whose proof is given in Section 5.

\section{$(-1)$-Curves}
In this section we classify all the $(-1)$-curves on the blow up of
$\p^2$ that may have negative ($\leq -2$) intersection with the elements
of a quasi-homogeneous linear system. With abuse of language we call
$(-1)$-curves also the curves on $\p^2$ whose proper transforms are
$(-1)$-curves. In order to find such curves we first prove that they must
be quasi-homogeneous. 
Then we restrict ourselves to the case of curves having points of
multiplicity at most 2 (since otherwise they would have intersection
$\geq -1$ with an element of a quasi-homogeneous linear system of
multiplicity 5), and we give a complete classification with the aid of
numerical properties.

\begin{propo}[\cite{cm1}]
\label{le:product}
Let $\ls$ be a linear system of $\p^2$ and let $E_1,E_2$ be two distinct 
$(-1)$-curves, such that $E_i\cdot\ls<0$ for $i\in\{1,2\}$. 
Then $E_1\cdot E_2=0$. 
\end{propo}
\begin{proof}
By hypothesis $E_1, E_2\in\bs\mid\ls\mid$.
Suppose that $E_1\cdot E_2\geq 1$. From the exact sequence
\[
0\rt H^0(S,E_2)\rt 
H^0(S,E_1+E_2)\rt
H^0(E_1,\oc_{E_1}(E_1
\cdot E_2-1))\rt 0,
\]
it follows that 
$h^0(S,E_1+E_2)>1$. This means that 
$E_1+E_2$ moves in a linear system, but this is impossible 
since $E_1+E_2\in\bs\mid\ls\mid$.
\end{proof}

Now, consider a curve $E=L(d,m_0,m_1,\ldots, m_r)$, 
with $m_i\geq 0$. 
Suppose that $E$ is a $(-1)$-curve such that $E\in\bs\mid\ls\mid$. 
Let $\sigma\in S_r$ be a permutation of the set $\{1,\ldots, r\}$ and consider
the curve $E_{\sigma}=L(d,m_0,m_{\sigma(1)},\ldots, m_{\sigma(r)})$.
Clearly $E_{\sigma}\cdot E_{\sigma}=E\cdot E=-1$ and $E_{\sigma}\cdot
K_{\p^2}=E\cdot K_{\p^2}=-1$, hence also $E_{\sigma}$ is a
$(-1)$-curve. Observe that since $\ls$ is quasi-homogeneous, then if 
$E\in\bs\mid\ls\mid$, by symmetry also $E_{\sigma}\in\bs\mid\ls\mid$. 
Now consider the permutation $\sigma_{ij}$ which switches $i$ with $j$ and 
leaves fixed all the other numbers. From Lemma \ref{le:product} it follows 
that $E\cdot E_{\sigma_{ij}}=0$. Hence we have:
\[
E^2=-1\hspace{2cm} E\cdot E_{\sigma_{ij}}=0.
\]
These equations translate into the system:
\[
\begin{cases}
d^2-m_0^2-\sum\limits_{k=1}^rm_k^2=-1 & \\
d^2-m_0^2-\sum\limits_{\stackrel{k=1}{k\neq i,j}}^rm_k^2-2m_im_j=0 & 
\end{cases}
\]
From these, it follows that $(m_i-m_j)^2=1$ which gives $m_i=m_j\pm 1$.
Hence, for each $i,j$, $m_i$ and $m_j$ are equal or they differ only by $1$. 
This means that the $(-1)$-curve $E$ may be rewritten as 
$L(d,m_0,m^s,(m+1)^{r-s})$. 
There are $\binom{r}{s}$ distinct and linearly independent 
curves of this kind. But $\dim\pic(S)=r+1$, so there are only four
possibilities: $s=0,\ s=1,\ s=r-1,\ s=r$. If $s=0$ or $s=r$ then the 
$(-1)$-curve $E$ belongs to a quasi-homogeneous system. 
We call the curve $E_{\rm tot}=\sum_{\sigma}E_{\sigma}$ 
{\em compound}. 

Hence we have proved the following: 

\begin{propo}
\label{propo:-1curves,1}
Let $\ls$ be a quasi-homogeneous linear system on $\p^2$ and let $E$ be a
$(-1)$-curve such that $E\cdot\ls<0$. Then $E$ is necessarily of the form:
\[
\begin{array}{ll}
E=L(d,m_0,m^r) &  \\
E=L(d,m_0,m-1,m^{r-1}) & E_{\rm tot}=L(rd,rm_0,(rm-1)^r) \\
E=L(d,m_0,m+1,m^{r-1}) & E_{\rm tot}=L(rd,rm_0,(rm+1)^r) \\
\end{array}
\]
\end{propo}

In order to classify $(-1)$-curves belonging to quasi-homogeneous systems, we
begin by classifying quasi-homogeneous curves of multiplicity $m\leq 2$. 
\begin{propo}[\cite{cm1}]
\label{propo:qo-1curves}
For each $m>1$, there are only a finite number of quasi-homogeneous 
$(-1)$-curves of multiplicity $m$. For $m\leq 2$, the quasi-homogeneous 
$(-1)$-curves are:
\[
L(1,1,1),\ L(2,0,1^5),\ L(e,e-1,1^{2e})\ L(6,3,2^7)
\]

\end{propo}
\begin{proof}
Let $E=L(d,m_0,m^r)$ be a $(-1)$-curve. Then $E$ must satisfy the following
equations:
\begin{equation}
\label{eq:-1curves,1}
\begin{cases}
d^2-m_0^2-rm^2=-1 & \\
-3d+m_0+rm=-1 & 
\end{cases}
\end{equation}
Eliminating $r$ from the two equations we obtain:
\begin{equation}
\label{eq:-1curves,2}
m-3dm+mm_0+d^2-m_0^2+1=0.
\end{equation}
If we put $x=2d-3m$ and $y=2m_0-m$, equation \eqref{eq:-1curves,2} becomes 
\begin{equation}
\label{eq:-1curves,3}
x^2-y^2=4(2m^2-m-1).
\end{equation}
This is the equation of an irreducible conic for $m\geq 2$ 
and it has a finite number of solutions. In fact,
the numbers $x-y$ and $x+y$ must be chosen among the divisors
of $4(2m^2-m-1)$ and there are only a finite number of choices.
In particular, for $m=2$ we obtain $x^2-y^2=20$, or equivalently $(d-3)^2-(m_0-1)^2=5$.
Hence we must have $d-m_0-2=1$ and $d+m_0-4=5$
which gives $d=6,\ m_0=3$ and the system $\ls(6,3,2^7)$.

Finally, for $m=1$, \eqref{eq:-1curves,2} becomes: $(d+m_0-2)(d-m_0-1)=0$.
This equation has the following solutions: 
$d=1,\ m_0=1$, which gives the system $\ls(1,1,1)$; $d=2, m_0=0$, 
which gives the system $\ls(2,0,1^5)$; $d=m_0+1$, corresponding to $\ls(d,d-1,1^{2d})$.
\end{proof}

\begin{propo}
The $(-1)$-curves of multiplicity $m\leq 2$ are listed in the following table:
\label{propo:-1curves,2}
\begin{table}[h]
\caption{$(-1)$-curves of multiplicity $\leq 2$}
\label{tab:-1planecurves}
\begin{tabular}{lc} 
\hline
$(-1)$-{\rm curves} $E$ & {\rm compound curve} $E_{\rm tot}$\\
\hline
  $L(2,0,1^5)$ &  \\ 
  $L(e,e-1,1^{2e})$ &  \\ 
  $L(6,3,2^7)$ &  \\ 
  $L(1,1,1)$ & $L(e,e,1^e)$  \\ 
  $L(1,0,1^2)$ & $L(3,0,2^3)$  \\ \hline
\end{tabular}
\end{table}
\end{propo}  
\begin{proof}
If $E=L(d,m_0,m-1,m^{r-1})$ then $E_{\rm tot}=L(rd,rm_0,(rm-1)^r)$ 
and it must be $rm-1\leq 2$. Since $m\geq 1$ (otherwise $m-1<0$)
and $r\geq 2$ (otherwise the system is not compound), the only possibilities 
are $m=1$ and $r=2,3$. 
In the first case, we obtain the $(-1)$-curve $E=L(1,1,1)$ which leads to the
system $E_{\rm tot}=L(r,r,1^r)$.
In the second case, $E=L(1,0,1^2)$ and $E_{\rm tot}=L(3,0,2^3)$. \\

If $E=L(d,m_0,m+1,m^{r-1})$ then $E_{\rm tot}=L(rd,rm_0,(rm+1)^r)$
and it must be $rm+1\leq 2$.
If $m=0$ then we obtain the $(-1)$-curve $E=L(1,1,1)$ which gives the compound 
curve: $E_{\rm tot}=L(r,r,1^r)$. 
If $m=1$ then $r\leq 1$ and there are not compound systems with
only two points. 
\end{proof}

\newpage

\section{Quasi-homogeneous systems of multiplicity $5$}

In this section we classify all the $(-1)$-special 
quasi-homogeneous linear systems $\ls$ of multiplicity $5$. 
\begin{propo}\label{-1special}
All the $(-1)$-special quasi-homogeneous systems of multiplicity $5$ are 
listed in the following table:

\begin{center}
\begin{table}[h]
\caption{$(-1)$-special system of multiplicity $5$}
\label{tab:specialsystems=4}
\begin{tabular}{lcrr}
\hline
System & & virtual dim. & effective dim. \\ \hline   
$\ls(d,d,5^r)$ & & $d-15r$ & $d-5r$ \\
$\ls(d,d-1,5^r)$ & & $2d-15r$ & $2d-9r$ \\
$\ls(d,d-2,5^r)$ & $(d,r)\neq(8e,2e)$ & $3d-15r-1$ & $3d-12r-1$\\
$\ls(8e,8e-2,5^{2e})$ & & $-6e-1$ & $0$ \\
$\ls(d,d-3,5^r)$ & $(d,r)\neq (7e,2e),$ & $4d-15r-3$ & $4d-14r-3$ \\
                 & $(7e+1,2e)$ & & \\
$\ls(7e,7e-3,5^{2e})$ & & $-2e-3$ & $0$ \\
$\ls(7e+1,7e-2,5^{2e})$ & & $-2e+1$ & $2$ \\
$\ls(8,0,5^2)$ & & $14$ & $15$ \\
$\ls(10,2,5^4)$ & & $2$ & $3$ \\
$\ls(3e+5,3e-2,5^{2e})$ & $1\leq e\leq 3$ & $19-6e$ & $20-6e$ \\
$\ls(3e+4,3e-3,5^{2e})$ & $1\leq e\leq 2$ & $11-6e$ & $14-6e$ \\
$\ls(4e+4,4e-2,5^{2e})$ & $1\leq e\leq 7$ & $13-2e$ & $14-2e$ \\
$\ls(4e+3,4e-3,5^{2e})$ & $1\leq e\leq 4$ & $6-2e$  & $9-2e$ \\
$\ls(4e+2,4e-4,5^{2e})$ & $1\leq e\leq 2$ & $-2e-1$ & $5-2e$ \\
$\ls(5e+3,5e-2,5^{2e})$ &           & $8$     & $9$ \\
$\ls(5e+2,5e-3,5^{2e})$ &           & $2$     & $5$ \\
$\ls(5e+1,5e-4,5^{2e})$ &           & $-4$    & $2$ \\
$\ls(5e,5e-5,5^{2e})$   &           & $-10$   & $0$ \\
$\ls(6e+2,6e-2,5^{2e})$ &           & $4$     & $5$ \\
$\ls(6e+1,6e-3,5^{2e})$ &           & $-1$    & $2$ \\
$\ls(6e,6e-4,5^{2e})$   &           & $-6$    & $0$ \\
$\ls(8,0,5^3)$          &           & $-1$    & $2$ \\
$\ls(8,1,5^3)$          &           & $-2$    & $1$ \\
$\ls(11,0,5^5)$         &           & $2$     & $5$ \\
$\ls(11,1,5^5)$         &           & $1$     & $4$ \\
$\ls(11,2,5^5)$         &           & $-1$     & $2$ \\ \hline
\end{tabular}
\end{table}
\end{center}
\end{propo}

\begin{proof}
The procedure for finding these $(-1)$-special systems may be described as follows: 

First step: find all the $(-1)$-curves $E_i$ such that 
\begin{equation}
\label{eq:specialsys,1}
E_i\cdot\ls=\mu_i\leq -2.
\end{equation}

Second step: let $\ms=\ls-\sum\mu_iE_i$ be the residual system. We need to 
verify that
\begin{equation}
\label{eq:specialsys,2}
v(\ms)\geq 0.
\end{equation}
Observe that to complete the first step, once we have found a $(-1)$-curve $E_i$
that satisfies \eqref{eq:specialsys,1}, we have to see if there is a
$(-1)$-curve $E_j$ such that $E_j\cdot (\ls-\mu_i E_i)\leq -2$ and so on.
We may speed the procedure by making an induction on the multiplicity of the
system $\ls$. Since the system $\ls-\mu_i E_i$ has multiplicity less than $4$,
we need only to know if this system is again $(-1)$-special or not. In the
affirmative case we search for the next $(-1)$-curve $E_j$, otherwise 
we proceed with step two.
So we need a complete list of quasi-homogeneous $(-1)$-special systems with $m
\leq 3$. The following table may be found in \cite{cm1}.
\begin{table}[h]
\caption{$(-1)$-special systems of multiplicity $\leq 3$}
\label{tab:specialsystems<4}
\begin{tabular}{lccc}
\hline
System & & virtual dim. & effective dim. \\
\hline
$\ls(4,0,2^5)$ & & $v=-1$ & $l = 0$ \\
$\ls(2e,2e-2,2^{2e})$, & $e \geq 1$ & $v = -1$ & $l = 0$\\
$\ls(d,d,2^e)$, & $d \geq 2e\geq 2$ & $v = d-3e$ & $l = d-2e$ \\
$\ls(4,0,3^2)$ & & $v=2$ & $l = 3$ \\
$\ls(6,0,3^5)$ & & $v=-3$ & $l=0$ \\
$\ls(6,2,3^4)$ & & $v=0$ & $l=1$ \\
$\ls(3e,3e-3,3^{2e})$, & $e\geq 1$ & $v = -3$ & $l=0$ \\
$\ls(3e+1,3e-2,3^{2e})$, & $e\geq 1$ &  $v=1$ & $l=2$ \\
$\ls(4e,4e-2,3^{2e})$, & $e\geq 1$ & $v = -1$ & $l=0$ \\
$\ls(d,d-1,3^e)$, &  $2d\geq 5e \geq 5$ & $v=2d-6e$ & $l = 2d-5e$ \\
$\ls(d,d,3^e)$, & $d \geq 3e \geq 3$ & $v = d-6e$ & $l = d-3e$ . \\ \hline
\end{tabular}
\end{table}

In order to satisfy \eqref{eq:specialsys,1} for a system $\ls=\ls(d,m_0,4^r)$,
the $(-1)$-curves $E_i$ may have multiplicity at most $2$.
Hence the $E_i$ are those listed in table \ref{tab:-1planecurves}.
So there are five possibilities to analyze. \\

$\bullet\  \ls(d,m_0,5^r)\cdot L(1,1,1)=-\mu,$ \\

In this case $m_0=d-5+\mu$ and the residual system is

\[ 
\ms=\ls(d-\mu r,d-\mu r-5+\mu,(5-\mu)^r).
\]

If $\mu =5$ then the system is $\ls(d,d,5^r)$. The residual system $\ms$ 
is non-special of dimension $d-5r$.\\

If $\mu =4$ then the system is $\ls(d,d-1,5^r)$. The residual system $\ms$ 
is non-special of dimension $2d-9r$.\\

If $\mu =3$ then the system is $\ls(d,d-2,5^r)$. The residual system $\ms$ 
may be $(-1)$-special only if $d-3r=2e$ and $r=2e$. 
This implies that the system is $\ls(8e,8e-2,5^{2e})$. In this case,
$\ms=\ls(2e,2e-2,2^{2e})$ is $(-1)$-special of dimension $0$. If $\ms$ is
non-special, then it has dimension $3d-12r-1$. \\

If $\mu =2$ then the system is $\ls(d,d-3,5^r)$. The residual system $\ms$ 
may be $(-1)$-special only if $r=2e$ and $d-2r=3e$, or $d-2r=3e+1$. 
In the first case one obtain the system $\ls(7e,7e-3,5^{2e})$ of
effective dimension $0$. In the second case one obtain the system 
$\ls(7e+1,7e-2,5^{2e})$ of effective dimension $2$.
If $\ms$ is non-special, then it has dimension $4d-14r-3$. \\

$\bullet\ \ls(d,m_0,5^{2e})\cdot L(e,e-1,1^{2e})=-\mu$ \\
 
This gives:
\begin{equation}
\label{eq:specialsys:poss.2,2}
e(d-m_0-10)+m_0+\mu=0
\end{equation}

The residual system is $\ms=\ls(d-\mu e,m_0-\mu (e-1),(5-\mu)^{2e})$. This
implies that $m_0-\mu (e-1)\leq d-\mu e$, which gives $\mu+m_0\leq d$.
On the other hand, from equation \eqref{eq:specialsys:poss.2,2} one deduce that
$d-m_0-10<0$, since otherwise the sum cannot be $0$.
So we have the inequality $\mu\leq d-m_0\leq + 9$. Let $k=d-m_0$, then
$\mu\leq k\leq 9$. Substituting in \eqref{eq:specialsys:poss.2,2} we obtain:
\begin{equation}
\label{eq:specialsys:poss.2,3}
m_0=(10-k)e-\mu,\hspace{2cm} d=(10-k)e-\mu+k.
\end{equation}
The residual system is:

\[
\ms=\ls((10-k-\mu)e-\mu+k,(10-k-\mu)e,(5-\mu)^{2e})
\]

If $k=9$ then $10-k-\mu <0$, so there are no systems. \\

If $k=8$ then $10-k-\mu\geq 0$ only if $\mu =2$. In this case 
$\ms = \ls(6,0,3^{2e})$ is non-special of dimension $27-12e$.  
This gives the systems $\ls(8,0,5^2)$ and $\ls(10,2,5^4)$.\\

If $k=7$ then $\mu\leq 3$. 
If $\mu=2$, the system $\ls(3e+5,3e-2,5^{2e})$ has virtual dimension $19-6e$. 
The residual system $\ms=\ls(e+5,e,3^{2e})$ is non-special of dimension $20-6e$. 
If $\mu=3$ the system $\ls(3e+4,3e-3,5^{2e})$ has virtual dimension $11-6e$.
The residual system $\ms=\ls(4,0,2^{2e})$ has dimension $14-6e$. \\

If $k=6$ then $\mu\leq 4$. 
If $\mu=2$ the system $\ls(4e+4,4e-2,5^{2e})$ has virtual dimension $13-2e$ 
and the residual system $\ms=\ls(2e+4,2e,3^{2e})$ is non-special of dimension 
$14-2e$.
If $\mu=3$ then the system $\ls(4e+3,4e-3,5^{2e})$ has virtual dimension
$6-2e$ and the residual system $\ms=\ls(e+3,e,2^{2e})$ is non-special of
dimension $9-2e$.
If $\mu=4$ then the system $\ls(4e+2,4e-4,5^{2e})$ has virtual dimension
$-2e-1$ and the residual system $\ms=\ls(2,0,1^{2e})$ is non special of
dimension $5-2e$.
\\

If $k=5$ then 
if $\mu=2$ the system $\ls(5e+3,5e-2,5^{2e})$ has virtual dimension $8$ 
and the residual system $\ms=\ls(3e+3,3e,3^{2e})$ is non-special of
dimension $9$. 
If $\mu=3$ the system $\ls(5e+2,5e-3,5^{2e})$ has virtual dimension $2$
and the residual system $\ms=\ls(2e+2,2e,2^{2e})$ is non-special of
dimension $5$.
If $\mu=4$ the system $\ls(5e+1,5e-4,5^{2e})$ has virtual dimension $-4$
and the residual system $\ms=\ls(e+1,e,1^{2e})$ is non-special of
dimension $2$.
If $\mu=5$ the system $\ls(5e,5e-5,5^{2e})$ has virtual dimension $-10$
and the residual system $\ms=\ls(0,0,0)$ is non-special of dimension $0$.
\\

If $k=4$ then 
if $\mu=2$ the system $\ls(6e+2,6e-2,5^{2e})$ has virtual dimension $4$
and the residual system $\ms=\ls(4e+2,4e,3^{2e})$ is non-special of
dimension $5$.
If $\mu=3$ the system $\ls(6e+1,6e-3,5^{2e})$ has virtual dimension $-1$
and the residual system $\ms=\ls(3e+1,3e,2^{2e})$ is non-special of
dimension $2$.
If $\mu=4$ the system $\ls(6e,6e-4,5^{2e})$ has virtual dimension $-6$
and the residual system $\ms=\ls(2e,2e,1^{2e})$ is non-special of
dimension $0$.
If $\mu=5$ the system $\ls(6e-1,6e-5,5^{2e})$ has virtual dimension $-11$
and the residual system is empty.\\

$\bullet\  \ls(d,m_0,5^7)\cdot\ls(6,3,2^7) =-2,$ \\

So $6d-3m_0-68=0$. There are no integer values of $d$ and $m_0$ that satisfy
this equation. \\

$\bullet\  \ls(d,m_0,5^5)\cdot L(2,0,1^5)=-\mu,$ \\

with $\mu\in\{2,3,4\}$. This equation gives $2d-25=-\mu$, so $\mu$ must be 
odd. If $\mu=3$ then $d=11$ and the residual system
$\ms=\ls(5,m_0,2^5)$ is not $(-1)$-special. Moreover $v(\ms)=5-m_0(m_0+1)/2$, so
$v(\ms)\geq 0$ only if $m_0\leq 2$. Hence the systems are: $\ls(11,m_0,5^5),\
m_0\in\{0,1,2\}$.\\

$\bullet\  \ls(d,m_0,5^3)\cdot L(1,0,1^2)=-2,$ \\

so $d=8$ and the residual system $\ms=\ls(2,m_0,1^3)$ has virtual dimension
$2-m_0(m_0+1)/2$. Hence it must be $m_0\leq 1$ and the systems are:
$\ls(8,0,5^3)$ and $\ls(8,1,5^3)$.
\end{proof}

\section{Degeneration of linear systems on $\p^2$}
In this section we recall the degeneration technique needed to specialize linear
systems on $\p^2$. For a reference see \cite{cm1}, Sections 2 and 3.\\

The idea is to take a flat family $X$ on a complex disk $\Delta$, such that
the fiber $X_t$ over a point $t\neq 0$ is a plane, while the central fiber $X_0$ 
is the union of a plane $\p$ and a Hirzebruch surface $\f=\f_1$. 
These two surfaces are 
joined transversely along a curve $R$ which is a line $L$ in $\p$ and the 
exceptional divisor $E$ on $\f$.

To give a linear system on $X_0$ is equivalent to giving two linear systems, on $\p$ and
on $\f$, which agree on the curve $R$.

Fix a positive integer $r$ and another non-negative integer $b\leq r$.
Let us consider $r-b+1$ general points $p_0, p_1, \ldots, p_{r-b}$ in $\p$
and $b$ general points $p_{r-b+1},...,p_r$ in $\f$.
These points are limits of $r+1$ general points  $p_{0,t}, p_{1,t},
\ldots, p_{r,t}$ in $X_t$. 
Let us call ${\ls}_t$ the linear system $\ls(d,m_0,m^r)$ in $X_t \cong
{\p}^2$ based at the points $p_{0,t}, p_{1,t},...,p_{r,t}$. 

For any integer $k$ there exists a linear system $\ls_0$ on $X_0$ that
restricts to $\p$ as a system $\lp=\ls(d-k,m_0,m^{r-b})$ and to $\f$ as a
system $\lf=\ls(d,d-k,m^b)$. One can prove that ${\ls_0}$ (for any $k$ and
$b$) can be obtained as a flat limit on $X_0$ of the system
${\ls}=\ls(d,m_0,m^r)$.

We say that ${\ls_0}$ is obtained from $\ls$
by a \emph{$(k,b)$-degeneration}.

We denote by $l_0$ the dimension of the linear system $\ls_0$ on $X_0$.
By semi-continuity, $l_0$
is not smaller than the dimension of the linear system on the general fiber,
i.e.,
\[
l_0 = \dim(\ls_0) \geq l = \dim \ls(d,m_0,m^r).
\]
Therefore we have the following:

\begin{lemma}[Lemma 2.1, \cite{cm1}]
\label{lo=Ethenreg}
If $l_0$ is equal to the expected dimension $v$
of $\ls$,
then the system $\ls$ is non-special.
\end{lemma}

Let $\klf$ and $\klp$ be the kernels of the restriction of the systems $\lf$
and $\lp$ to $R$, while $\rs_{\f}$ and $\rs_{\p}$ are the restricted systems.
The dimension $l_0$ is obtained in terms of the dimensions of the systems 
$\lp$ and $\lf$, and the dimensions of the subsystems $\klp \subset \lp$
and $\klf \subset \lf$ consisting of divisors containing the double curve $R$.
Notice that by slightly abusing notation we have
\[
\begin{array}{ll}
\lp = \ls(d-k,m_0,m^{r-b}) & \klp = \ls(d-k-1,m_0,m^{r-b}) \\ 
\lf = \ls(d,d-k,m^b) & \klf = \ls(d,d-k+1,m^b)
\end{array}
\]

We recall the following notations:
\[
\begin{array}{ll}
v = v(d,m_0,m^r) &
\text{ the virtual dimension of $\ls$,} \\
\vp = v(d-k,m_0,m^{r-b}) &
\text{ the virtual dimension of $\lp$ }, \\
\vf = v(d,d-k,m^b) &
\text{ the virtual dimension of $\lf$ }, \\
\kvp = v(d-k-1,m_0,m^{r-b}) &
\text{ the virtual dimension of $\klp$}, \\
\kvf = v(d,d-k+1,m^b) &
\text{ the virtual dimension of $\klf$ }, \\
l = l(d,m_0,m^r) &
\text{ the dimension of $\ls$}, \\
\dlp = l(d-k,m_0,m^{r-b}) &
\text{ the dimension of $\lp$ }, \\
\dlf = l(d,d-k,m^b) &
\text{ the dimension of $\lf$ }, \\
\dklp = l(d-k-1,m_0,m^{r-b}) &
\text{ the dimension of $\klp$ }, \\
\dklf = l(d,d-k+1,m^b) &
\text{ the dimension of $\klf$ }, \\
r_{\p} = \dlp-\dklp-1 &
\text{ the dimension of the restricted} \\ 
& \text{  system $\rs_{\p}$}\\
r_{\f} = \dlf-\dklf-1 &
\text{ the dimension of the restricted} \\
& \text{ system $\rs_{\f}$}\\
\end{array}
\]

Then, with the help of a transversality lemma of linear systems on $R$,
one can prove the following Proposition on the dimension $l_0$ that we are
going to use in next section.

\begin{propo}[Proposition 3.3, \cite{cm1}]
\label{pro:dimLo}
The following identities hold:
\begin{itemize}
\item[(a)] If $r_{\p}+r_{\f} \leq d-k-1$, then $l_0 = \dklp + \dklf + 1$.
\item[(b)] If $r_{\p}+r_{\f} \geq d-k-1$, then $l_0 = \dlp+ \dlf - d+k$.
\item[(c)] $v = \vp+\vf -d+k = \vf + \kvp +1 = \vp + \kvf +1$
\end{itemize}
\end{propo}

\section{main theorem}

In this section we prove the main theorem on quasi-homogeneous linear
systems of multiplicity 5, using the degeneration techniques and 
induction on the number of points.

\begin{mainthm}
\label{teo:main}
$\ls(d,m_0,5^r)$ is special if and only if it is $(-1)$-special.
\end{mainthm}

Before proving the main theorem, we state some technical lemmas that
simplify the proof of the theorem. The proof of these lemmas is given in 
next section.

\begin{lemma}
\label{le:threepoints}
A linear system with at most three base points is special if and only
if it is $(-1)$-special.
\end{lemma}

\begin{lemma}
\label{le:m_0>=d-7}
If $m_0\geq d-7$, then the system $\ls(d,m_0,5^r)$ is special if and only
if it is $(-1)$-special.
\end{lemma}

\begin{lemma}
\label{le:d<=100}
If $d\leq 150$, then the system $\ls(d,m_0,5^r)$ is special if and only if
it
is $(-1)$-special.
\end{lemma}

\begin{proof}[Proof of Main Theorem]
By Lemma \ref{le:m_0>=d-7}, we may suppose that $d\geq m_0+8$.
We proceed by induction on $r$.\\
Perform a $(4,b)$-degeneration obtaining the systems:
\[
\begin{array}{ll}
\lp=\ls(d-4,m_0,5^{r-b}) & \lf=\ls(d,d-4,5^b) \\
\klp=\ls(d-5,m_0,5^{r-b}) & \klf=\ls(d,d-3,5^b) 
\end{array}  
\]
Step I:\ $v(\ls(d,m_0,5^r))\leq -1$\\

\begin{claim}\label{claim1}
If we can find a $(4,b)$-degeneration such that 
\begin{equation}\label{5cose}
\dklf=-1,\ \kvp\leq v,\ \dklp=\kvp,\ \vf=\dlf,\ \vp=\dlp
\end{equation}
then the system $\ls(d,m_0,5^r)$ is empty.
\end{claim}
In fact we have that $\dklf=-1$ and $\dklp=\kvp\leq v\leq -1$, and hence
$\klf=\klp=\emptyset$. Moreover 
\[
\begin{split}
r_{\p}+r_{\f}&=\dlp-\dklp-1+\dlf-\dklf-1\\
&=\dlp+\dlf\\
&=\vp+\vf\\
&=v+d-4\leq d-5, 
\end{split}
\]
By Proposition
\ref{pro:dimLo} and Lemma \ref{lo=Ethenreg} we have that the system is empty
and hence the claim follows.
\vskip .2truecm
We find conditions on the integer $b$ necessary to guarantee 
\eqref{5cose}.

First of all, $\klf=\ls(d,d-3,5^b)$ is a $(-1)$-special system, and it is empty if
$4d-14b-3\leq -1$, i.e. $b\geq \frac{2d-1}{7}$.

Moreover, $\kvp-v=15b-5d+5\leq 0$ if and only if $b\leq\frac{d-1}{3}$.

Looking at the table \ref{tab:specialsystems=4}
we have that the system  $\lf=\ls(d,d-4,5^{b})$ can be 
special only if it is one of the following systems
\[
\begin{split}
&\ls(6e+2,6e-2,5^{2e}),\\
&\ls(6e+1,6e-3,5^{2e}),\\
&\ls(6e,6e-4,5^{2e}).
\end{split}
\] 
But if $b_0$ is an integer such that $\ls(d,d-4,5^{b_0})=\ls(6e+\epsilon,
6e+\epsilon-4,5^{2e})$, then taking $b_0+\lambda$, with $\lambda\geq 1$, 
we obtain a non-special system.

Concerning the system $\lp=\ls(d-4,m_0,5^{r-b})$, 
since we are supposing that $m_0\leq d-8$ (and then 
the difference between the degree and $m_0$ is at least 4),
it can not be special if $d\geq 12$ and  
$r-b$ is odd. 

For the same reason we can say that $\klp=\ls(d-5,m_0,5^{r-b})$ is non-special
if $d\geq 12$ and $r-b$ is odd. The only remaining case is $m_0=d-8$, because
$\klp=\ls(d-5,d-8,5^{r-b})$ is $(-1)$-special, and its dimension is 
$\dklp=4(d-5)-14(r-b)-3$. In this case, instead of proving that $\dklp=\kvp$, 
we look for a $b$ such that $\dklp=4(d-5)-14(r-b)-3<0$. This is equivalent
to saying that $b<r+\frac{23}{14}-\frac{2d}{7}$. Moreover we know that the 
linear system we started with was $\ls(d,d-8,5^r)$, which has negative 
dimension only if $r>\frac{3d}{5}-\frac{28}{15}$. Then, if
$b<\frac{3d}{5}-\frac{28}{15}+\frac{23}{14}-\frac{2d}{7}=\frac{66d-47}{210}$,
the linear system $\klp$ is empty.
\vskip .2truecm
Summarizing, if the difference 
\[ 
\frac{66d-47}{210}-\frac{2d-1}{7}
\]
is at least 4, i.e. $d\geq 143$, we can choose two values $b=b_0$ or $b_0+2$ 
between $\frac{2d-1}{7}$
and $\frac{66d-47}{210}$, and such that $r-b$ is odd. Clearly 
for one of them the system $\ls(d,d-4,5^{b})$ is non-special.

Hence all the conditions of \eqref{5cose} hold, and the theorem is true
in the case $v\leq -1$.
\vskip .2truecm
Step II:\ $v(\ls(d,m_0,5^r))\geq 0$

In this situation the following holds.
\begin{claim}
\label{claim2}
If $\lp$ and $\lf$ are non-special and $v-1\geq\dklp+\dklf$, then 
$v=l_0$.
\end{claim}
In fact, we have that 
\begin{equation}
\begin{split}
r_{\p}+r_{\f}-(d-k-1)&=\dlp-\dklp-1+\dlf-\dklf-1-(d-k-1)\\
&=\vp+\vf-d+k-1-(\dklp+\dklf)\\
&=v-1-(\dklp+\dklf)\geq 0,
\end{split}
\end{equation}
which means that $r_{\p}+r_{\f}\geq d-k-1$ and, by Proposition \ref{pro:dimLo},
$l_0=\dlp+\dlf-d+k$. But since we are supposing that $\lp$ and $\lf$
are non-special, $l_0=\vp+\vf-d+k=v$, and hence the claim.
\vskip .2truecm
As before, if $r-b$ is odd and $d\geq 12$, then $\lp$ is non-special. Moreover,
if $b$ can be choosen in a sufficiently large interval, also $\lf$ is
non-special.

The system $\klf$ is $(-1)$-special and $\dklf-\kvf=b$. If $m_0<d-8$, $r-b$ 
is odd and $d\geq 12$, then $\klp$ is non-special. This implies that
$\dklp+\dklf=\kvp+\kvf+b$. We have to guarantee that $\kvp+\kvf+b\leq v-1
=\vf+\kvp$, which is equivalent to $b\leq\vf-\kvf=d-3$.

In the case $m_0=d-8$, $\klp=\ls(d-5,d-8,5^{r-b})$ is $(-1)$-special and
$\dklp=\kvp+r-b$, which implies $\dklp+\dklf=\kvp+\kvf+r$. If we want this 
sum to be smaller than or equal to $v-1=\kvp+\vf$, it must be $r\leq\vf-
\kvf=d-3$. But since the starting system was $\ls=\ls(d,d-8,5^r)$ and 
we are supposing that $v(\ls)\geq 0$, we have that $r\leq \frac{3d}{5}
-\frac{28}{15}$. In particular, if $d\geq\frac{17}{6}$, then $r\leq d-3$,
and we are in the hypothesis of Claim \ref{claim2}.

Summarizing, if $d\geq 12$, we can find a $b$ such that 
the hypothesis of Claim \ref{claim2} hold, and hence the theorem in the case
$v(\ls)\geq 0$.
\end{proof}

\section{Proofs of lemmas}

\begin{proof}[Proof of Lemma \ref{le:threepoints}]
In order to prove this, it is sufficient to consider the three points
$(1:0:0),(0:1:0),(0:0:1)$ and then evaluate the condition imposed by these
points on the monomials $x_0^{a_0}x_1^{a_1}x_2^{a_2}$ where
$a_0+a_1+a_2=$degree.
\end{proof}

\begin{proof}[Proof of Lemma \ref{le:m_0>=d-7}]

The case $m_0\geq d-m-1$, for quasi-homogeneous systems of multiplicity
$m$, has been proved in \cite{cm1}. Hence we have only to prove the case $m_0=d-7$. \\

Performing recursively $k$ Cremona transformations, we obtain 
the system $\ls(d-3k,d-3k-7,5^{r-2k},2^{2k})$. Let us put
$d-7=3t+\epsilon$, with  $\epsilon=0,1,2$, and $r=2q+\eta$, with $\eta=0,1$.
\vskip .2truecm
\noindent $\bullet\  t\leq q$.
We perform $t$ transformations to obtain
the system $\ls(7+\epsilon,\epsilon,5^{r-2t},2^{2t})$.
\begin{enumerate}
\parindent =0pt
\item $\epsilon=0$

The system is $\ls(7,5^{r-2t},2^{2t})$. If $r-2t\geq 2$, then the line 
passing through two of the points with multiplicity $5$ is a
$(-1)$-curve contained three times in the system. Then either
$\ls(7,5^{r-2t},2^{2t})$ is $(-1)$-special, or it is empty. 

If $r-2t=0,1$, then the linear system is homogeneous or 
quasi-homogeneous of multiplicity $2$, and hence the conjecture
holds.
\item $\epsilon=1$

As before, if $r-2t\geq 2$, the linear system $\ls(8,1,5^{r-2t},2^{2t})$
is $(-1)$-special or empty.

If $r-2t=0,1$, the linear system is equivalent to a homogeneous
or a quasi-homogeneous linear system of multiplicity $2$, and the
conjecture holds.

\item $\epsilon=2$

We obtain the linear system $\ls(9,5^{r-2t},2^{2t+1})$.

If $r-2t\geq 4$, the system is $(-1)$-special or empty.

If $r-2t=3$, if we cut away the three lines passing through the
points of multiplicity $5$, we obtain the equivalent 
system $\ls(6,3^3,2^{2t})$, and then, after another Cremona 
transformation, $\ls(3,2^{2t})$, for which the conjecture holds.

If $r-2t=2$, we cut away the line passing through the points
of multiplicity $5$, and then we perform a Cremona transformation to
obtain $\ls(6,2^{2t+1})$. The conjecture still holds.

\end{enumerate}
\vskip .2truecm
$\bullet\  t>q$.
We can perform $q$ Cremona transformations and get $\ls(d-3q,
d-7-3q,5^{\eta},2^{2q})$.

If $\eta=0$, the linear system is quasi-homogeneous of
multiplicity $2$, and the conjecture holds.

If $\eta=1$ the linear system is $\ls(\delta,\delta-7,5,2^{2q})$, with $\delta
=d-3q$. 
Perform a $(1,b)$-degeneration obtaining the systems:
\[
\begin{array}{ll}
\lp=\ls(\delta-1,\delta-7,5,2^{2q-b}) & \lf=\ls(\delta,\delta-1,2^b) \\
\klp=\ls(\delta-2,\delta-7,5,2^{2q-b}) & \klf=\ls(\delta,\delta,2^b) 
\end{array}  
\]
We proceed exactly as we did in the proof of Main Theorem. 
\vskip .2truecm
Step I:\ $v(\ls(d,d-7,5^{2q+1}))\leq -1$\\

We want to find a $b$ such that 
\begin{equation}\label{5cose+}
\dklf=-1,\ \kvp\leq v,\ \dklp=\kvp,\ \vf=\dlf,\ \vp=\dlp.
\end{equation}
Let us find out the conditions on the integer $b$ necessary to guarantee 
\eqref{5cose+}.

Looking at the table \ref{tab:specialsystems=4} we have that the system $\lf$ 
is non-special (because it is not $(-1)$-special), while $\klf$ is
$(-1)$-special, of dimension $\dklf=\delta-2b$. In particular,
if $b\geq \frac{\delta+1}{2}$, $\klf$ is empty.

Concerning $\lp=\ls(\delta-1,\delta-7,5,2^{2q-b})$ and 
$\klp=\ls(\delta-2,\delta-7,5,2^{2q-b})$, we can perform some Cremona 
transformations and obtain that they are not special. Moreover, 
$\kvp-v=6b-4\delta-2$ is not bigger than zero if we choose $b\leq
\frac{2\delta+1}{3}$.

Summarizing, if we can choose an integer $b$ such that
\[
\frac{\delta+1}{2}\leq b\leq\frac{2\delta+1}{3},
\]
i.e. if the difference $\frac{2\delta+1}{3}-\frac{\delta+1}{2}$ is at least 1,
then we can guarantee all the hypothesis of \eqref{5cose+}. But this is 
equivalent to saying that $\delta\geq 7$, or $d-7\geq 3q$, which is true since 
$d-7=3t+\epsilon$, and we are supposing $t>q$.
\vskip .2truecm

Step II:\ $v(\ls(d,d-7,5^{2q+1}))\geq 0$\\

In this case, Claim \ref{claim2} still holds, and hence, since
we have already seen that $\lp$ and $\lf$ are non-special, to end 
the proof of the lemma we only have to see that 

\begin{equation}
\label{v-1}
v-1\geq\dklp+\dklf.
\end{equation}

The system $\klf$ is $(-1)$-special and $\dklf-\kvf=b$, while $\klp$ is 
non-special. Summing up, $\dklp+\dklf=\kvp+\kvf+b$, while $v-1=\kvp+\vf$ and 
hence \eqref{v-1} is equivalent to $b\leq\vf-\kvf=\delta=d-3q$. 

Therefore, if we take an integer $b\leq d-3q$, and perform a 
$(1,b)$-degeneration, we are in the hypothesis of the claim.
\end{proof}

\begin{proof}[Proof of Lemma \ref{le:d<=100}]

We prove this lemma by writing an algorithm that performs 
degeneration on quasi-homogeneous linear systems.
First we define a function that evaluates the effective dimension 
of a linear system $\ls$ according to the conjecture. 
We do this in a recursive way: if $\ls$ has negative product with some $(-1)$-curve $E$
(i.e. $\ls\cdot E=-t<0$), then we redefine $\ls$ as $\ls-tE$ and we restart from the beginning.
Otherwise we define the effective dimension of $\ls$ as its virtual dimension.
If $\ls$ is empty, it may happen that the procedure never ends, since at each step
there is a $(-1)$-curve $E$ such that $\ls\cdot E<0$. In this case after a finite 
number of steps the system $\ls$ has a negative degree.
This may happen only if the system is empty, in which case we set the 
effective dimension of $\ls$ to be $-1$. \\

Then we define a test function, that tries all the possible degenerations 
of a linear system $\ls$ and returns the value $0$ if there exists a degeneration such that
the dimension of $\ls_0$ is equal to the virtual dimension of $\ls$.
The evaluation of $\dim\ls_0$ is performed with the aid of Proposition \ref{pro:dimLo}.
Due to the induction hypothesis the dimension of $\lp,\ \lf,\ \klp,\ \klf$ may be 
evaluated with the recursive function already defined. \\

In the final step we consider all the quasi-homogeneous linear systems 
$\ls(d,m_0,5^r)$ such that $8\leq d\leq 150$ and $m_0\leq d-8$. For each 
one of them all the possible degenerations are performed with the test function. \\

This procedure was written in \textsc{Maple} and may be found at \\
\verb1http://socrates.mat.unimi.it/~laface1

We give here a list of the systems for which the degeneration does not work:

\begin{center}
\begin{table}[ht]
\caption{The exceptional cases}
\label{tab:exc}
\begin{tabular}{clllll}
\hline
Degree & System & & & &\\ \hline   
12 & $\ls(12,\alpha,5^5),$ & $\ls(12,\beta,5^6);$    & $1\leq\alpha\leq 4$ & $0\leq\beta\leq 4$ \\
13 & $\ls(13,\alpha,5^6),$ & $\ls(13,\beta,5^7);$    & $\alpha=4,5$        & $\beta=0,1$ \\
14 & $\ls(14,\alpha,5^7),$ & $\ls(14,\beta,5^8);$    & $\alpha=4,5,6$      & $\beta=0,1$ \\ 
   & $\ls(14,6,5^6),$      &                         &                     & \\ 
16 & $\ls(16,\alpha,5^9),$ & $\ls(16,\beta,5^{10});$ & $\alpha=5,6$        & $\beta=0,1,2$ \\ 
17 & $\ls(17,8,5^9);$      &                         &                     & \\
18 & $\ls(18,10,5^9);$     &                         &                     & \\
19 & $\ls(19,9,5^{11}),$   & $\ls(19,11,5^{10});$    &                     & \\
20 & $\ls(20,11,5^{10}),$  & $\ls(20,12,5^{10});$    &                     & \\
21 & $\ls(21,13,5^{11});$  &                         &                     & \\
23 & $\ls(23,14,5^{13});$  &                         &                     & \\
26 & $\ls(26,18,5^{14});$  &                         &                     & \\
28 & $\ls(28,20,5^{15});$  &                         &                     & \\
33 & $\ls(33,25,5^{18}).$  &                         &                     & \\
\hline
\end{tabular}
\end{table}
\end{center}

Almost all these systems can be studied with quadratic transformations. 
But for some of them it has been necessary to evaluate the Hilbert function 
of the corresponding ideal.
We did this evaluation with the aid of \textsc{Singular} \cite{sing1}.
The source of our program can be downloaded at the following url:
\verb1http://socrates.mat.unimi.it/~laface1

Even in these cases the conjecture was true.
\end{proof}

\bibliographystyle{amsalpha}

\providecommand{\bysame}{\leavevmode\hbox to3em{\hrulefill}\thinspace}

\end{document}